\documentclass[reqno]{amsart}

\newtheorem{theorem}{Theorem}
\newtheorem{lemma}[theorem]{Lemma}
\newtheorem{proposition}[theorem]{Proposition}
\newtheorem{corollary}[theorem]{Corollary}

\theoremstyle{definition}
\newtheorem{definition}[theorem]{Definition}
\usepackage{graphicx}
\usepackage[all]{xy}
\usepackage{pstricks,enumerate,pst-node,pst-text,pst-plot,stmaryrd}

\theoremstyle{remark}

\numberwithin{equation}{section}
\numberwithin{theorem}{section}

\newcommand{\intav}[1]{\mathchoice {\mathop{\vrule width 6pt height 3 pt depth  -2.5pt
\kern -8pt \intop}\nolimits_{\kern -6pt#1}} {\mathop{\vrule width
5pt height 3  pt depth -2.6pt \kern -6pt \intop}\nolimits_{#1}}
{\mathop{\vrule width 5pt height 3 pt depth -2.6pt \kern -6pt
\intop}\nolimits_{#1}} {\mathop{\vrule width 5pt height 3 pt depth
-2.6pt \kern -6pt \intop}\nolimits_{#1}}}

\newcommand{\intavl}[1]{\mathchoice {\mathop{\vrule width 6pt height 3 pt depth  -2.5pt
\kern -8pt \intop}\limits_{\kern -6pt#1}} {\mathop{\vrule width 5pt
height 3  pt depth -2.6pt \kern -6pt \intop}\nolimits_{#1}}
{\mathop{\vrule width 5pt height 3 pt depth -2.6pt \kern -6pt
\intop}\nolimits_{#1}} {\mathop{\vrule width 5pt height 3 pt depth
-2.6pt \kern -6pt \intop}\nolimits_{#1}}}

\newcommand{\mc}{\mathcal}

\newcommand{\B}{\mc{B}}
\newcommand{\cc}{\mc{C}}

\newcommand{\F}{\mc{F}}

\newcommand{\M}{\mc{M}}

\newcommand{\orb}{\mc{O}}

\newcommand{\R}{\mathbb{R}}
\newcommand{\N}{\mathbb{N}}

\newcommand{\Z}{\mathbb{Z}}
\newcommand{\X}{\mathbb{X}}

\newcommand{\Folner}{F$\phi$lner }
\newcommand{\cyl}{\text{\rm Cyl}}

\newcommand{\fr}{\text{\rm fr}}

\newcommand{\al}{\alpha}
\newcommand{\be}{\beta}

\begin{document}

\title[$\Z^d$-actions with prescribed topological and ergodic properties]
{$\Z^d$-actions with prescribed topological and ergodic properties}

\author[Yuri Lima]{Yuri Lima}
\address{Instituto Nacional de Matem\'atica Pura e Aplicada, Estrada Dona Castorina 110, 22460-320, Rio de Janeiro, Brasil.}
\email{yurilima@impa.br}


\subjclass[37B05, 37B10, 37B40]{37B05, 37B10, 37B40}

\date{October, 15, 2010}


\begin{abstract}
We extend constructions of Hahn-Katznelson \cite{HK} and Pavlov \cite{P} to $\Z^d$-actions on symbolic
dynamical spaces with prescribed topological and ergodic properties. More specifically, we describe a
method to build $\Z^d$-actions which are (totally) minimal, (totally) strictly ergodic and have
positive topological entropy.
\end{abstract}

\maketitle

\section{Introduction}

Ergodic theory studies statistical and recurrence properties of measurable transformations $T$ acting in a
probability space $(X,\B,\mu)$, where $\mu$ is a measure invariant by $T$, that is, $\mu(T^{-1}A)=\mu(A)$,
for all $A\in\B$. It investigates a wide class of notions, such as ergodicity, mixing and entropy. These
properties, in some way, give qualitative and quantitative aspects of the randomness of $T$. For example,
ergodicity means that $T$ is indecomposable in the metric sense with respect to $\mu$ and entropy
is a concept that counts the exponential growth rate for the number of statistically significant
distinguishable orbit segments.

In most cases, the object of study has topological structures: $X$ is a compact metric space,
$\B$ is the Borel $\sigma$-algebra of $X$, $\mu$ is a Borel measure probability and $T$ is a homeomorphism of
$X$. In this case, concepts such as minimality and topological mixing give topological aspects of the randomness
of $T$. For example, minimality means that $T$ is indecomposable in the topological sense, that is, the orbit of every
point is dense in $X$.

A natural question arises: how do ergodic and topological concepts relate to each other? How do ergodic properties
forbid topological phenomena and vice-versa? Are metric and topological indecomposability equivalent? This last
question was answered negatively in \cite{F} via the construction of a minimal diffeomorphism of the torus
$\mathbb T^2$ which preserves area but is not ergodic.

Another question was raised by W. Parry: suppose $T$ has a unique Borel probability invariant
measure and that $(X,T)$ is a minimal transformation. Can $(X,T)$ have positive entropy? The difficulty in
answering this at the time was the scarcity of a wide class of minimal and uniquely ergodic transformations.
This was solved affirmatively in \cite{HK},  where F. Hahn and Y. Katznelson developed an inductive method
of constructing symbolic dynamical systems with the required properties. The principal idea of the paper was
the weak law of large numbers.

Later, works of Jewett and Krieger (see \cite{Pe}) proved that every ergodic measure-preserving system
$(X,\B,\mu, T)$ is metrically isomorphic to a minimal and uniquely ergodic homeomorphism on a Cantor set
and this gives many examples to Parry's question: if an ergodic system $(X,\B,\mu, T)$ has positive metric
entropy and $\Phi:(X,\B,\mu,T)\rightarrow (Y,\mathcal C,\nu,S)$ is the metric isomorphism obtained by
Jewett-Krieger's theorem, then $(Y,S)$ has positive topological entropy, by the variational principle.

It is worth mentioning that the situation is quite different in smooth ergodic theory, once some regularity on
the transformation is assumed. A. Katok showed in \cite{K} that every $C^{1+\al}$ diffeomorphism of a
compact surface can not be minimal and simultaneously have positive topological entropy. More specifically,
he proved that the topological entropy can be written in terms of the exponential growth of periodic points
of a fixed order.

Suppose that $T$ is a mesure-preserving transformation on the probability space $(X,\B,\mu)$ and
$f:X\rightarrow \R$ is a measurable function. A successful area in ergodic theory deals with the convergence
of averages $n^{-1}\cdot\sum_{k=1}^n f\left(T^kx\right)$, $x\in X$, when $n$ converges to infinity.
The well known {\it Birkhoff's Theorem} states that such limit exists for almost every $x\in X$ whenever
$f$ is an $L^1$-function. Several results have been (and still are being) proved when, instead of $\{1,2,\ldots,n\}$,
average is made along other sequences of natural numbers. A remarkable result on this direction was given by J.
Bourgain \cite{Bo}, where he proved that if $p(x)$ is a polynomial with integer coefficients and $f$ is an
$L^p$-function, for some $p>1$, then the averages $n^{-1}\cdot\sum_{k=1}^n f\left(T^{p(k)}x\right)$ converge for
almost every $x\in X$. In other words, convergence fails to hold for a negligible set with respect to the measure
$\mu$. In \cite{B}, V. Bergelson asked if this set is also negligible from the topological point of view. It turned
out, by a result of R. Pavlov \cite{P}, that this is not true. He proved that, for every sequence
$(p_n)_{n\ge 1}\subset\Z$ of zero upper-Banach density, there exist a totally minimal, totally uniquely
ergodic and topologically mixing transformation $(X,T)$ and a continuous function $f:X\rightarrow\R$ such that
$n^{-1}\cdot\sum_{k=1}^n f\left(T^{p_k}x\right)$ fails to converge for a residual set of $x\in X$.

Suppose now that $(X,T)$ is totally minimal, that is, $(X,T^n)$ is minimal for every positive integer $n$.
Pavlov also proved that, for every sequence $(p_n)_{n\ge 1}\subset\Z$ of zero upper-Banach density, there
exists a totally minimal, totally uniquely ergodic and topologically mixing continuous transformation $(X,T)$
such that $x\not\in\overline{\{T^{p_n}x\,;\,n\ge 1\}}$ for an uncountable number of $x\in X$.

In this work, we extend the results of Hahn-Katznelson and Pavlov, giving a method of constructing (totally)
minimal and (totally) uniquely ergodic $\Z^d$-actions with positive topological entropy. We carry out our program
by constructing closed shift invariant subsets of a sequence space. More specifically, we build a sequence of finite
configurations $(\cc_k)_{k\ge 1}$ of $\{0,1\}^{\Z^d}$, $\cc_{k+1}$ being essentially formed by the
concatenation of elements in $\cc_k$ such that each of them occurs
statistically well-behaved in each element of $\cc_{k+1}$, and consider the set of limits of shifted
$\cc_k$-configurations as $k\rightarrow+\infty$. The main results are

\begin{theorem}\label{thm 1}
There exist totally strictly ergodic $\Z^d$-actions $(X,\B,\mu,T)$ with arbitrarily large positive topological
entropy.
\end{theorem}

\noindent We should mention that this result is not new, because Jewett-Krieger's Theorem is true for
$\Z^d$-actions \cite{We}. This formulation emphasizes to the reader that the constructions,
which may be used in other settings, have the additional advantage of controlling the topological entropy.

\begin{theorem}\label{thm 3}
Given a set $P\subset\Z^d$ of zero upper-Banach density, there exist a totally strictly ergodic
$\Z^d$-action $(X,\B,\mu,T)$ and a continuous function
$f:X\rightarrow\R$ such that the ergodic averages
$$\dfrac{1}{|P\cap(-n,n)^d|}\sum_{g\in P\cap(-n,n)^d} f\left(T^gx\right)$$
fail to converge for a residual set of $x\in X$. In addition, $(X,\B,\mu,T)$ can have arbitrarily large
topological entropy.
\end{theorem}

The above theorem has a special interest when $P$ is an {\it arithmetic set} for which classical ergodic theory
and Fourier analysis have established almost-sure convergence. This is the case (also proved in \cite{Bo}) when
$$P=\{(p_1(n),\ldots,p_d(n))\,;\,n\in\Z\},$$
where $p_1,\ldots,p_d$ are polynomials with integer coefficients: for any $f\in L^p$, $p>1$, the limit
$$\lim_{n\rightarrow+\infty}\dfrac{1}{n}\sum_{k=1}^n f\left(T^{(p_1(k),\ldots,p_d(k))}x\right)$$
exists almost-surely. Note that $P$ has zero upper-Banach density whenever one
of the polynomials has degree greater than 1.

\begin{theorem}\label{thm 2}
Given a set $P\subset\Z^d$ of zero upper-Banach density, there exists a totally strictly ergodic
$\Z^d$-action $(X,\B,\mu,T)$ and an uncountable set $X_0\subset X$ for which
$x\not\in\overline{\{T^{p_n}x\,;\,n\ge 1\}}$, for every $x\in X_0$. In addition, $(X,\B,\mu,T)$
can have arbitrarily large topological entropy.
\end{theorem}

Yet in the arithmetic setup, Theorem \ref{thm 2} is the best topological result one can
expect. Indeed, Bergelson and Leibman proved in \cite{BL} that if $T$ is a minimal $\Z^d$-action, then there
is a residual set $Y\subset X$ for which $x\in\overline{\{T^{(p_1(n),\ldots,p_d(n))}x\,;\,n\in\Z\}}$, for every
$x\in Y$.

\section{Preliminaries}

We begin with some notation. Consider a metric space $X$, $\B$ its Borel $\sigma$-algebra and a $G$ group
with identity $e$. Throughout this work, $G$ will denote $\Z^d$, $d>1$, or one of its subgroups.

\subsection{Group actions}

\begin{definition}
A $G$-action on $X$ is a measurable transformation $T:G\times X\rightarrow X$,
denoted by $(X,T)$, such that
\begin{enumerate}[(i)]
\item $T(g_1,T(g_2,x))=T(g_1g_2,x)$, for every $g_1,g_2\in G$ and $x\in X$.
\item $T(e,x)=x$, for every $x\in X$.
\end{enumerate}
\end{definition}

\noindent In other words, for each $g\in G$, the restriction
$$
\begin{array}{ccrcl}
T^g&:&X&\longrightarrow &X\\
 & & x&\longmapsto &T(g,x)\\
\end{array}
$$
is a bimeasurable transformation on $X$ such that $T^{g_1g_2}=T^{g_1}T^{g_2}$, for
every $g_1,g_2\in G$. When $G$ is abelian, $(T^{g})_{g\in G}$ forms a commutative group of bimeasurable
transformations on $X$. For each $x\in X$, the {\it orbit} of $X$ with respect to $T$ is the set
$$\orb_T(x) \doteq\{T^gx\,;\,g\in G\}.$$
If $F$ is a subgroup of $G$, the restriction $T|_{F}:F\times X\rightarrow X$ is clearly
a $F$-action on $X$.

\begin{definition}
We say that $(X,T)$ is {\it minimal} if $\orb_T(x)$ is dense in $X$, for every $x\in X$, and
{\it totally minimal} if $\orb_{T|_F}(x)$ is dense in $X$, for every $x\in X$ and every
subgroup $F<G$ of finite index.
\end{definition}

\noindent Remind that the {\it index} of a subgroup $F$, denoted by $(G:F)$, is the number of cosets
of $F$ in $G$. The above definition extends the notion of total minimality
of $\Z$-actions. In fact, a $\Z$-action $(X,T)$ is totally minimal
if and only if $T^n:X\rightarrow X$ is a minimal transformation, for every $n\in\Z$.

Consider the set $\M(X)$ of all Borel probability measures in $X$. A probability $\mu\in\M(X)$ is
{\it invariant under} $T$ or simply $T$-{\it invariant} if
$$\mu\left(T^gA\right)=\mu(A),\ \forall\,g\in G,\ \forall\,A\in\B.$$
Let $\M_T(X)\subset\M(X)$ denote the set of all $T$-invariant probability measures. Such set is non-empty
whenever $G$ is amenable, by a Krylov-Bogolubov argument applied to any \Folner sequence of $G$.

\begin{definition}
A {\it $G$ measure-preserving system} or simply $G$-{\it mps} is a quadruple $(X,\B,\mu,T)$, where $T$ is a $G$-action
on $X$ and $\mu\in\M_T(X)$.
\end{definition}

We say that $A\in\mathcal B$ is $T$-invariant if $T^gA=A$, for all $g\in G$.

\begin{definition}
The $G$-mps $(X,\B,\mu,T)$ is {\it ergodic} if it has only trivial invariant sets, that is, if $\mu(A)=0$ or $1$ whenever
$A$ is a measurable set invariant under $T$.
\end{definition}

\begin{definition}
The $G$-action $(X,T)$ is {\it uniquely ergodic} if $\M_T(X)$ is unitary, and
{\it totally uniquely ergodic} if, for every subgroup $F<G$ of finite index, the restricted $F$-action
$(X,T|_F)$ is uniquely ergodic.
\end{definition}

\begin{definition}
We say that $(X,T)$ is {\it strictly ergodic} if it is minimal and uniquely ergodic, and {\it totally strictly
ergodic} if, for every subgroup $F<G$ of finite index, the restricted $F$-action
$(X,T|_F)$ is strictly ergodic.
\end{definition}

The result below was proved in \cite{W} and states the pointwise ergodic theorem for $\Z^d$-actions.

\begin{theorem}
Let $(X,\B,\mu,T)$ be a $\Z^d$-mps. Then, for every $f\in L^1(\mu)$, there is a $T$-invariant
function $\tilde f\in L^1(\mu)$ such that
$$\lim_{n\rightarrow+\infty}\dfrac{1}{n^d}\sum_{g\in[0,n)^d}f\left(T^gx\right)=\tilde f(x)$$
for $\mu$-almost every $x\in X$. In particular, if the action is ergodic, $\tilde f$ is constant and
equal to $\int fd\mu$.
\end{theorem}

Above, $[0,n)$ denotes the set $\{0,1,\ldots,n-1\}$, $[0,n)^d$ the $d$-dimensional cube
$[0,n)\times\cdots\times[0,n)$ of $\Z^d$ and by a $T$-invariant function we mean that
$f\circ T^g=f$, for every $g\in G$. These averages allow the characterization of unique
ergodicity. Let $C(X)$ denote the space of continuous functions from $X$ to $\R$.

\begin{proposition}\label{unique ergodicity}
Let $(X,T)$ be a $\Z^d$-action on the compact metric space $X$. The following items are equivalent.
\begin{enumerate}[(a)]
\item $(X,T)$ is uniquely ergodic.
\item For every $f\in C(X)$ and $x\in X$, the limit
$$\lim_{n\rightarrow+\infty}\dfrac{1}{n^d}\sum_{g\in[0,n)^d}f\left(T^gx\right)$$
exists and is independent of $x$.
\item For every $f\in C(X)$, the sequence of functions
$$f_n=\dfrac{1}{n^d}\sum_{g\in[0,n)^d}f\circ T^g$$
converges uniformly in $X$ to a constant function.
\end{enumerate}
\end{proposition}

\begin{proof}
The implications (c)$\Rightarrow$(b)$\Rightarrow$(a) are obvious. It remains to prove (a)$\Rightarrow$(c).
Let $\M_T(X)=\{\mu\}$. We'll show that $f_n$ converges uniformly to $\tilde f=\int fd\mu$. By contradiction, suppose
this is not the case for some $f\in C(X)$.
This means that there exist $\varepsilon >0$, $n_i\rightarrow\infty$ and $x_i\in X$ such that
$$\left|f_{n_i}(x_i)-\int fd\mu\right|\ge\varepsilon.$$
For each $i$, let $\nu_i\in\mathcal M(X)$ be the probability measure associated to the linear functional
$\Theta_i:C(X)\rightarrow\R$ defined by
$$\Theta_i(\varphi)=\frac{1}{{n_i}^d}\sum_{g\in[0,n_i)^d}\varphi(T^gx_i)\,,\ \varphi\in C(X).$$
Restricting to a subsequence, if necessary, we assume that $\nu_i\rightarrow\nu$ in the weak-star topology.
Because the cubes $A_i=[0,n_i)^d$ form a \Folner sequence in $\Z^d$, $\nu\in\mathcal M_T(X)$.
In fact, for each $h\in\Z^d$,
\begin{eqnarray*}
\left|\int\left(\varphi\circ T^h\right)d\nu-\int\varphi d\nu\right|
&=&\lim_{i\rightarrow\infty}\frac{1}{{n_i}^d}\left|\sum_{g\in A_i+h}\varphi(T^gx_i)-\sum_{g\in A_i}\varphi(T^gx_i)\right|\\
&\le&\max_{x\in X}{|\varphi(x)|}\cdot\lim_{i\rightarrow\infty}\dfrac{\#A_i\Delta(A_i+h)}{\#A_i}\\
&=&0.
\end{eqnarray*}
But
$$\left|\int fd\nu-\int fd\mu\right|=\lim_{i\rightarrow\infty}\left|\int fd\nu_i-\int fd\mu\right|
=\lim_{i\rightarrow\infty}\left|f_{n_i}(x_i)-\int fd\mu\right|\ge\varepsilon$$
and so $\nu\not=\mu$, contradicting the unique ergodicity of $(X,T)$.
\end{proof}

\subsection{Subgroups of $\Z^d$}\label{subsection subgroups}

Let $\F$ be the set of all subgroups of $\Z^d$ of finite index. This set is countable, because each element of $\F$ is generated by
$d$ linearly independent vectors of $\Z^d$. Consider, then, a subset $(F_k)_{k\in\N}$ of $\F$ such that, for each
$F\in\F$, there exists $k_0>0$ such that $F_k<F$, for every $k\ge k_0$. For this, just consider an enumeration
of $\F$ and define $F_k$ as the intersection of the first $k$ elements.
Such intersections belong to $\F$ because
$$(\Z^d:F\cap F')\le (\Z^d:F)\cdot (\Z^d:F'),\ \forall\, F,F'<\Z^d.$$
Restricting them, if necessary, we assume that $F_k=m_k\cdot\Z^d$,
where $(m_k)_{k\ge 1}$ is an increasing sequence of positive integers. Observe that
$(\Z^d:F_k)={m_k}^d$. Such sequence will be fixed throughout the rest of the paper.

\begin{definition} Given a subgroup $F<\Z^d$, we say that two elements
$g_1,g_2\in\Z^d$ are {\it congruent modulo} $F$ if $g_1-g_2\in F$ and denote it by $g_1\equiv_F g_2$.
The set $\bar F\subset\Z^d$ is a {\it complete residue set modulo $F$} if, for every $g\in\Z^d$, there
exists a unique $h\in\bar F$ such that $g\equiv_F h$.
\end{definition}

\noindent Every complete residue set modulo $F$ is canonically identified to the quocient $\Z^d/F$ and
has exactly $(\Z^d:F)$ elements.

\subsection{Symbolic spaces}\label{subsection simbolic space}

Let $\cc$ be a finite alphabet and consider the set $\Omega(\cc)={\cc}^{\Z^d}$
of all functions $x:\Z^d\rightarrow\cc$. We endow $\cc$ with the discrete topology and
$\Omega(\cc)$ with the product topology. By Tychonoff's theorem, $\Omega(\cc)$ is a compact metric
space. We are not interested in a particular metric in $\Omega(\cc)$. Instead, we consider a basis of
topology $\B_0$ to be defined below.

Consider the family $\mathcal R$ of all finite $d$-dimensional cubes $A=[r_1,r_1+n)\times\cdots\times[r_d,r_d+n)$
of $\Z^d$, $n\ge 0$. We say that $A$ has length $n$ and is centered at $g=(r_1,\ldots,r_d)\in\Z^d$.

\begin{definition}
A {\it configuration} or {\it pattern} is a pair $b_A=(A,b)$, where
$A\in\mathcal R$ and $b$ is a function from $A$ to $\cc$.
We say that $b_A$ is {\it supported in $A$ with encoding function $b$}.
\end{definition}

Let $\Omega_A(\cc)$ denote the space of configurations supported in $A$ and $\Omega^*(\cc)$
the space of all configurations in $\Z^d$:
$$\Omega^*(\cc)\doteq\{b_A\,;\, b_A\text{ is a configuration}\}.$$
Given $A\in\mathcal R$, consider the map $\Pi_A:\Omega(\mathcal C)\rightarrow\Omega_A(\mathcal C)$
defined by the restriction
$$
\begin{array}{ccrcl}
\Pi_A(x)&:&A&\longrightarrow &\mathcal C\\
 & & g&\longmapsto &x(g)\\
\end{array}
$$
In particular, $\Pi_{\{g\}}(x)=x(g)$. We use the simpler notation $x|_A$ to denote $\Pi_A(x)$.

\begin{definition}\label{defi 1}
If $A\in\mathcal R$ is centered at $g$, we say that $x|_A$ is a {\it configuration of $x$ centered at $g$}
or that $x|_A$ {\it occurs in $x$ centered at $g$}.
\end{definition}

For $A_1,A_2\in\mathcal R$ such that $A_1\subset A_2$, let $\pi_{A_1}^{A_2}:\Omega_{A_2}\rightarrow\Omega_{A_1}$
be the restriction
$$
\begin{array}{ccrcl}
\pi_{A_1}^{A_2}(b)&:&A_1&\longrightarrow & \mathcal C\\
 & & g&\longmapsto &b(g)\\
\end{array}
$$
As above, when there is no ambiguity, we denote $\pi_{A_1}^{A_2}(b)$ simply by $b|_{A_1}$.
It is clear that the diagram below commutes.
$$
\xymatrix{
\Omega(\mathcal C)\ar[r]^{\Pi_{A_2}}\ar[dr]_{\Pi_{A_1}}&\Omega_{A_2}(\mathcal C)\ar[d]^{\pi_{A_1}^{A_2}}\\
&\Omega_{A_1}(\mathcal C)}
$$
These maps will help us to control the patterns to appear in the constructions of Section
\ref{section constructions}.

By a cylinder in $\Omega(\mathcal C)$ we mean the set of elements of $\Omega(\mathcal C)$ with some fixed
configuration. More specifically, given $b_A\in\Omega^*(\mathcal C)$, the {\it cylinder generated by} $b_A$
is the set
$$\cyl(b_A)\doteq\{x\in\Omega(\cc)\,;\,x|_A=b_A\}.$$
The family $\B_0:=\{\cyl(b_A)\,|\,b_A\in\Omega^*(\mathcal C)\}$
forms a clopen set of cylinders generating $\B$. Hence the set
$C_0=\{\chi_B\,;\,B\in\B_0\}$ of cylinder characteristic functions generates a dense subspace in
$C(\Omega(\cc))$. Let $\mu$ be the probability measure defined by
$$\mu(\cyl(b_A))=|\cc|^{-|A|},\ \forall\,b_A\in\Omega^*(\mathcal C)\,,$$
and extended to $\B$ by Carathe\'odory's Theorem. Above, $|\cdot|$ denotes the number of elements of a set.

Consider the $\Z^d$-action $T:\Z^d\times\Omega(\mathcal C)\rightarrow\Omega(\mathcal C)$ defined by
$$T^g(x)=(x(g+h))_{h\in\Z^d}\,,$$
also called the {\it shift action}. Given $B=\cyl(b_A)$ and $g\in\Z^d$, let $B+g$ denote the cylinder
associated to $b_{A+g}=(\tilde b,A+g)$, where $\tilde b:A+g\rightarrow\{0,1\}$ is defined by
$\tilde b(h)=b(h-g)$, $\forall\,h\in A+g$. With this notation,
\begin{equation}\label{equacao 1}
\chi_B\circ T^g=\chi_{B+g}\,.
\end{equation}
In fact,
\begin{eqnarray*}
\chi_B(T^gx)&=&1\\
\iff\hspace{2cm}T^gx&\in& B\\
\iff\hspace{1.3cm}x(g+h)&=& b(h),\ \forall\,h\in A\\
\iff\hspace{1.9cm}x(h)&=&\tilde b(h),\ \forall\,h\in A+g\\
\iff\hspace{2.4cm}x&\in&B+g.
\end{eqnarray*}

\begin{definition}
A {\it subshift} of $(\Omega(\mathcal C),T)$ is a $\Z^d$-action $(X,T)$, where $X$ is a closed
subset of $\Omega(\mathcal C)$ invariant under $T$.
\end{definition}

\subsection{Topological entropy}

For each subset $X$ of $\Omega(\mathcal C)$ and $A\in\mathcal R$, let
$$\Omega_A(\mathcal C,X)=\{x|_A\,;\,x\in X\}$$
denote the set of configurations supported in $A$
which occur in elements of $X$ and $\Omega^*(\mathcal C,X)$ the space of all configurations in
$\Z^d$ occuring in elements of $X$,
$$\Omega^*(\mathcal C,X)=\bigcup_{A\in\mathcal R}\Omega_A(\mathcal C,X).$$

\begin{definition}
The {\it topological entropy} of the subshift $(X,T)$ is the limit
\begin{equation}\label{entropy eq 1}
h(X,T)=\lim_{n\rightarrow+\infty}\dfrac{\log{|\Omega_{[0,n)^d}(\mathcal C,X)}|}{n^d}\,,
\end{equation}
which always exists and is equal to $\inf_{n\in\N}\frac{1}{n^d}\cdot\log{|\Omega_{[0,n)^d}(\mathcal C,X)}|$.
\end{definition}

\subsection{Frequencies and unique ergodicity}

\begin{definition}
Given configurations $b_{A_1}\in\Omega_{A_1}(\mathcal C)$ and $b_{A_2}\in\Omega_{A_2}(\mathcal C)$,
the {\it set of ocurrences of $b_{A_1}$ in $b_{A_2}$} is
$$S(b_{A_1},b_{A_2})\doteq\{g\in\Z^d\,;\,A_1+g\subset A_2\text{ and }\pi_{A_1+g}^{A_2}(b_{A_2})=b_{A_1+g}\}.$$
The {\it frequency of $b_{A_1}$ in $b_{A_2}$} is defined as
$$\fr(b_{A_1},b_{A_2})\doteq\dfrac{|S(b_{A_1},b_{A_2})|}{|A_2|}\,\cdot$$
Given $F\in\F$ and $h\in\Z^d$, the {\it set of ocurrences of $b_{A_1}$ in $b_{A_2}$ centered at
$h$ modulo $F$} is
$$S(b_{A_1},b_{A_2},h,F)\doteq\{g\in S(b_{A_1},b_{A_2})\,;\,A_1+g\text{ is centered at a vertex }\equiv_F h\}$$
and the {\it frequency of $b_{A_1}$ in $b_{A_2}$ centered at $h$ modulo $F$} is the quocient
$$\fr(b_{A_1},b_{A_2},h,F)\doteq\dfrac{|S(b_{A_1},b_{A_2},h,F)|}{|A_2|}\,\cdot$$
\end{definition}

\noindent Observe that if $\bar F\subset\Z^d$ is a complete residue set modulo $F$, then
$$\fr(b_{A_1},b_{A_2})=\sum_{g\in\bar F}\fr(b_{A_1},b_{A_2},g,F).$$

To our purposes, we rewrite Proposition \ref{unique ergodicity} in a different manner.

\begin{proposition}
A subshift $(X,T)$ is uniquely ergodic if and only if, for every $b_A\in\Omega^*(\cc)$ and $x\in X$,
$$\fr(b_A,x)\doteq\lim_{n\rightarrow+\infty}\fr\left(b_A,x|_{[0,n)^d}\right)$$
exists and is independent of $x$.
\end{proposition}

\begin{proof}
By approximation, condition (b) of Proposition \ref{unique ergodicity} holds for $C(X)$ if and only
if it holds for $C_0=\{\chi_B\,;\,B\in\B_0\}$. If $f=\chi_{\cyl(b_A)}$, (\ref{equacao 1}) implies that
\begin{eqnarray*}
\lim_{n\rightarrow+\infty}f_n(x)&=&\lim_{n\rightarrow+\infty}\dfrac{1}{n^d}\sum_{g\in[0,n)^d}f\left(T^gx\right)\\
      &=&\lim_{n\rightarrow+\infty}\dfrac{1}{n^d}\sum_{g\in[0,n)^d}\chi_{\cyl\left(b_{A+g}\right)}(x)\\
      &=&\lim_{n\rightarrow+\infty}\dfrac{1}{n^d}\sum_{g\in[0,n)^d\atop{A+g\subset[0,n)^d}}\chi_{\cyl\left(b_{A+g}\right)}(x)\\
			&=&\lim_{n\rightarrow+\infty}\fr\left(b_A,x|_{[0,n)^d}\right)\\
			&=&\fr(b_A,x),
\end{eqnarray*}
where in the third equality we used that, for a fixed $A\in\mathcal R$,
$$\lim_{n\rightarrow+\infty}\dfrac{|\{g\in[0,n)^d\,;\,A+g\not\subset[0,n)^d\}|}{n^d}=0.$$
\end{proof}

\begin{corollary}
A subshift $(X,T)$ is totally uniquely ergodic if and only if, for every $b_A\in\Omega^*(\cc)$,
$x\in X$ and $F\in\F$,
$$\fr(b_A,x,F)\doteq\lim_{n\rightarrow+\infty}\fr\left(b_A,x|_{[0,n)^d},0,F\right)$$
exists and is independent of $x$.
\end{corollary}

So, unique ergodicity is all about constant frequencies. We'll obtain this via the Law of
Large Numbers, equidistributing ocurrences of configurations along residue classes of subgroups.

\subsection{Law of Large Numbers}

Intuitively, if $A$ is a subset of $\Z^d$, each letter of $\cc$ appears in $x|_A$ with frequency approximately
$1/|\cc|$, for almost every $x\in\Omega(\cc)$. This is what the Law of Large Number says. For our purposes, we state
this result in a slightly different way. Let $(X,\B,\mu)$ be a probability space and $A\subset\Z^d$ infinite.
For each $g\in A$, let $\X_g:X\rightarrow\R$ be a random variable.

\begin{theorem}\label{Law of large numbers}
(Law of Large Numbers) If $(\X_g)_{g\in A}$ is a family of independent and identically distributed
random variables such that $\mathbb E[\X_g]=m$, for every $g\in A$, then the sequence
$\left(\overline{\X}_n\right)_{n\ge 1}$ defined by
$$\overline{\X}_n=\dfrac{\sum_{g\in A\cap[0,n)^d}\X_g}{|A\cap[0,n)^d|}$$
converges in probability to $m$, that is, for any $\varepsilon>0$,
$$\lim_{n\rightarrow+\infty}\mu\left(\left|\overline{\X}_n-m\right|<\varepsilon\right)=1.$$
\end{theorem}

Consider the probability measure space $(X,\B,\mu)$ defined in Subsection
\ref{subsection simbolic space}. Fixed $w\in\cc$, let $\X_g:\Omega(\cc)\rightarrow\R$ be defined as
\begin{equation}\label{random variable 1}
\begin{array}{ccll}
\X_g(x)&=&1\ ,&\text{if }x(g)=w\\
&=&0\ ,&\text{if }x(g)\not=w.\\
\end{array}
\end{equation}
It is clear that $(\X_g)_{g\in\Z^d}$ are independent, identically distributed and satisfy
$$\mathbb E[\X_g]=\int_X\X_g(x)d\mu(x)=\dfrac{1}{|\cc|}\ ,\ \forall\,g\in\Z^d.$$
In addition,
$$\overline{\X}_n(x)=\dfrac{\sum_{g\in[0,n)^d}\X_g(x)}{n^d}=\dfrac{\left|S\left(w,x|_{[0,n)^d}\right)\right|}{n^d}=
\fr\left(w,x|_{[0,n)^d}\right),$$
which implies the

\begin{corollary}\label{law of large numbers}
Let $w\in\cc$, $g\in\Z^d$, $F\in\F$ and $\varepsilon>0$.
\begin{enumerate}
\item[(a)] The number of elements $b\in\Omega_{[0,n)^d}(\cc)$ such that
$$\left|\fr(w,b)-\dfrac{1}{|\cc|}\right|<\varepsilon$$
is assymptotic to $|\cc|^{n^d}$ as $n\rightarrow+\infty$.
\item[(b)] The number of elements $b\in\Omega_{[0,n)^d}(\cc)$ such that
$$\left|\fr(w,b,g,F)-\dfrac{1}{|\cc|\cdot(\Z^d:F)}\right|<\varepsilon$$
is assymptotic to $|\cc|^{n^d}$ as $n\rightarrow+\infty$.
\item[(c)] The number of elements $b\in\Omega_{[0,n)^d}(\cc)$ such that
$$\left|\fr(w,b,g,F)-\dfrac{1}{|\cc|\cdot(\Z^d:F)}\right|<\varepsilon$$
for every $w\in\cc$ and $g\in\Z^d$ is assymptotic to $|\cc|^{n^d}$ as $n\rightarrow+\infty$.
\end{enumerate}
\end{corollary}

\begin{proof}
(a) The required number is equal to
$$|\cc|^{n^d}\cdot\mu\left(\{x\in\Omega(\cc)\,;\,
\left|\fr\left(w,x|_{[0,n)^d}\right)-|\cc|^{-1}\right|<\varepsilon\}\right)$$
and is asymptotic to $|\cc|^{n^d}$, as the above $\mu$-measure converges to $1$.

\noindent (b) Take $A=F+g$ and $(\X_h)_{h\in A}$ as in (\ref{random variable 1}). For any $x\in\Omega(\cc)$,
\begin{eqnarray*}
\overline{\X}_n(x)&=&\dfrac{\left|S\left(w,x|_{[0,n)^d},g,F\right)\right|}{\left|A\cap[0,n)^d\right|}\\
                 &=&\fr\left(w,x|_{[0,n)^d},g,F\right)\cdot(\Z^d:F)+o(1)\,,
\end{eqnarray*}
because $\left|A\cap[0,n)^d\right|$ is asymptotic to $n^d/(\Z^d:F)$. This implies that for $n$ large
$$\left|\fr\left(w,x|_{[0,n)^d},g,F\right)-\dfrac{1}{|\cc|\cdot(\Z^d:F)}\right|<\varepsilon\ \iff\
\left|\overline{\X}_n(x)-\dfrac{1}{|\cc|}\right|<\varepsilon\cdot(\Z^d:F)$$
and then Theorem \ref{Law of large numbers} guarantees the conclusion.\\
\noindent (c) As the events are independent, this follows from (b).
\end{proof}

\section{Main Constructions}\label{section constructions}

Let $\cc=\{0,1\}$. In this section, we construct subshifts $(X,T)$ with
topological and ergodic prescribed properties. To this matter, we build a sequence
of finite non-empty sets of configurations $\cc_k\subset\Omega_{A_k}(\cc)$, $k\ge 1$, such that:
\begin{enumerate}[(i)]
\item $A_k=[0,n_k)^d$, where $(n_k)_{k\ge 1}$ is an increasing sequence of positive integers.
\item $n_1=1$ and $\cc_1=\Omega_{A_1}(\cc)\cong\{0,1\}$.
\item $\cc_k$ is the concatenation of elements of $\cc_{k-1}$, possibly with the insertion of few additional
blocks of zeroes and ones.
\end{enumerate}
Given such sequence $(\cc_k)_{k\ge 1}$, we consider $X\subset\Omega(\cc)$ as the set of limits
of shifted $\cc_k$-patterns as $k\rightarrow+\infty$, that is, $x\in X$ if there exist sequences
$(w_k)_{k\ge 1}$, $w_k\in\cc_k$, and $(g_k)_{k\ge 1}\subset\Z^d$ such that
$$x=\lim_{k\rightarrow+\infty}T^{g_k}w_k.$$
The above limit has an abuse of notation, because $T$ acts in $\Omega(\cc)$ and $w_k\not\in\Omega(\cc)$.
Formally speaking, this means that, for each $g\in\Z^d$, there exists $k_0\ge 1$ such that
$$x(g)=w_k(g+g_k),\ \forall\,k\ge k_0.$$
By definition, $X$ is invariant under $T$ and, for any $k$, every $x\in X$ is an
infinite concatenation of elements of $\cc_k$ and additional blocks of zeroes and ones.

If $\cc_k\subset\Omega_{A_k}(\{0,1\})$ and $A\in\mathcal R$, $\Omega_A(\cc_k)$ is identified in a natural
way to a subset of $\Omega_{n_kA}(\{0,1\})$. In some situations, to distinguish this association,
we use small letters for $\Omega_A(\cc_k)$ and capital letters for
$\Omega_{n_kA}(\{0,1\})$\footnote{For example, $w\in\Omega_A(\cc_k)$ and
$W\in\Omega_{n_kA}(\{0,1\})$ denote the ``same'' element.}. In this situation, if $w\in\Omega_A(\cc_k)$
and $g\in A$, the pattern $w(g)\in\cc_k$ occurs in $W\in\Omega_{n_kA}(\{0,1\})$ centered at $n_kg$.
In other words, if $w_k\in\cc_k$, then
\begin{equation}\label{equacao 9}
S(w_k,W,n_kg,F)=n_k\cdot S(w_k,w,g,F).
\end{equation}

In each of the next subsections, $(\cc_k)_{k\ge 1}$ is constructed with specific combinatorial and statistical
properties.

\subsection{Minimality}

The action $(X,T)$ is minimal if and only if, for each $x,y\in X$, every configuration
of $x$ is also a configuration of $y$. For this,
suppose $\cc_k\subset\Omega_{A_k}(\{0,1\})$ is defined and non-empty.

By the Law of Large Numbers, if $l_k$ is large,
every element of $\cc_k$ occurs in almost every element of $\Omega_{[0,l_k)^d}(\cc_k)$ (in fact, by
Corollary \ref{law of large numbers}, each of them occurs approximately with frequency $1/|\cc_k|>0$).
Take any subset $\cc_{k+1}$ of $\Omega_{[0,l_k)^d}(\cc_k)$ with this property and consider it as a subset
of $\Omega_{[0,n_{k+1})^d}(\{0,1\})$, where $n_{k+1}=l_kn_k$.

Let us prove that $(X,T)$ is minimal. Consider $x,y\in X$ and $x|_A$ a finite configuration
of $x$. For large $k$, $x|_A$ is a subconfiguration of some $w_k\in\cc_k$. As $y$ is formed by the
concatenation of elements of $\cc_{k+1}$, every element of $\cc_k$ is a configuration of $y$. In particular,
$w_k$ (and then $x|_A$) is a configuration of $y$.

\subsection{Total minimality}
The action $(X,T)$ is totally minimal if and only if, for each $x,y\in X$ and $F\in\F$, every
configuration $x|_A$ of $x$ centered\footnote{Because of the $T$-invariance of $X$, we can suppose
that $x|_A$ is centered in $0\in\Z^d$. In fact, instead of $x,y$, we consider $T^gx,T^gy$.} at $0$ also
occurs in $y$ centered at some $g\in F$. To guarantee this for every $F\in\F$, we
inductively control the ocurrence of subconfigurations centered in finitely many subgroups of $\Z^d$.

Consider the sequence $(F_k)\subset\F$ defined in Subsection \ref{subsection subgroups}.
By induction, suppose $\cc_k\subset\Omega_{A_k}(\{0,1\})$ is non-empty satisfying (i), (ii), (iii)
and the additional assumption
\begin{enumerate}
\item[(iv)] ${\rm gcd}(n_k,m_k)=1$ (observe that this holds for $k=1$).
\end{enumerate}
Take $l_k$ large and $\tilde\cc_{k+1}\subset\Omega_{[0,l_km_{k+1})^d}(\cc_k)$ non-empty such that
\begin{enumerate}
\item[(v)] $S(w_k,w|_{[0,l_km_{k+1}-1)^d},g,F_k)\not=\emptyset$, for every triple
$(w_k,w,g)\in \cc_k\times\tilde\cc_{k+1}\times\Z^d$.
\end{enumerate}
Considering $w|_{[0,l_km_{k+1}-1)^d}$ as an element of $\Omega_{[0,l_km_{k+1}n_k-n_k)^d}(\{0,1\})$,
(\ref{equacao 9}) implies that
$$S(w_k,W|_{[0,l_km_{k+1}n_k-n_k)^d},n_kg,F_k)\not=\emptyset,\,\forall\,(w_k,w,g)\in \cc_k\times \tilde\cc_{k+1}\times\Z^d.$$
As ${\rm gcd}(n_k,m_k)=1$, the set $n_k\Z^d$ runs over all residue classes modulo $F_k$ and so
(the restriction to $[0,l_km_{k+1}n_k-n_k)^d$ of) every element of $\tilde\cc_{k+1}$ contains every
element of $\cc_k$ centered at every residue class modulo $F_k$.

Obviously, $\cc_{k+1}$ must not be equal to $\tilde\cc_{k+1}$, because $m_{k+1}$ divides $l_km_{k+1}n_k$.
Instead, we take $n_{k+1}=l_km_{k+1}n_k+1$ and insert positions $B_i$, $i=1,2,\ldots,d$, next
to faces of the cube $[0,l_km_{k+1}n_k)^d$. These are given by
$$B_i=\{(r_1,\ldots,r_d)\in A_{k+1}\,;\,r_i=l_km_{k+1}n_k-n_k\}.$$
There is a natural surjection $\Phi:\Omega_{A_{k+1}}(\{0,1\})\rightarrow\Omega_{[0,n_{k+1}-1)^d}(\{0,1\})$
obtained removing the positions $B_1,\ldots,B_d$. More specifically, if
\begin{eqnarray*}
\delta(r)&=&0\,,\ \text{ if }r<l_km_{k+1}n_k-n_k\,,\\
         &=&1\,,\ \text{ otherwise}
\end{eqnarray*}
and 
\begin{equation}\label{equacao 10}
\Delta(r_1,\ldots,r_d)=(\delta(r_1),\ldots,\delta(r_d)),
\end{equation}
the map $\Phi$ is given by
$$\Phi(W)(g)=W\left(g+\Delta(g)\right),\ \forall\,(r_1,\ldots,r_d)\in[0,n_{k+1}-1)^d.$$
We conclude the induction step taking $\cc_{k+1}=\Phi^{-1}(\tilde\cc_{k+1})$.\\

\begin{center}
\psset{unit=.4cm} \begin{pspicture}(-.2,-2.8)(29,10.2)

\pspolygon(0,0)(9,0)(9,9)(0,9)\pspolygon(17,0)(27,0)(27,10)(17,10)

\psline[linewidth=.5pt](0,3)(9,3)\psline[linewidth=.5pt](0,6)(9,6)
\psline[linewidth=.5pt](3,0)(3,9)\psline[linewidth=.5pt](6,0)(6,9)
\psline[linewidth=.5pt](17,3)(23,3)\psline[linewidth=.5pt](17,6)(23,6)
\psline[linewidth=.5pt](23,0)(23,6)\psline[linewidth=.5pt](20,0)(20,6)
\pspolygon[linewidth=.5pt](24,0)(24,6)(27,6)(27,7)(24,7)(24,10)(23,10)(23,7)(17,7)(17,6)(23,6)(23,0)(24,0)
\psline[linewidth=.5pt](24,3)(27,3)\psline[linewidth=.5pt](20,7)(20,10)

\uput[0](1.4,-.8){$\Phi(w_{k+1})\in\tilde\cc_{k+1}$}\uput[0](19.1,-.8){$w_{k+1}\in\cc_{k+1}$}
\psline[linewidth=.3pt]{<-}(11,4.5)(15,4.5)
\uput[0](12.3,5.2){$\Phi$}
\uput[0](1,1.5){1}\uput[0](4,1.5){2}\uput[0](7,1.5){3}
\uput[0](1,4.5){4}\uput[0](4,4.5){5}\uput[0](7,4.5){6}
\uput[0](1,7.5){7}\uput[0](4,7.5){8}\uput[0](7,7.5){9}

\uput[0](18,1.5){1}\uput[0](21,1.5){2}\uput[0](25,1.5){3}
\uput[0](18,4.5){4}\uput[0](21,4.5){5}\uput[0](25,4.5){6}
\uput[0](18,8.5){7}\uput[0](21,8.5){8}\uput[0](25,8.5){9}
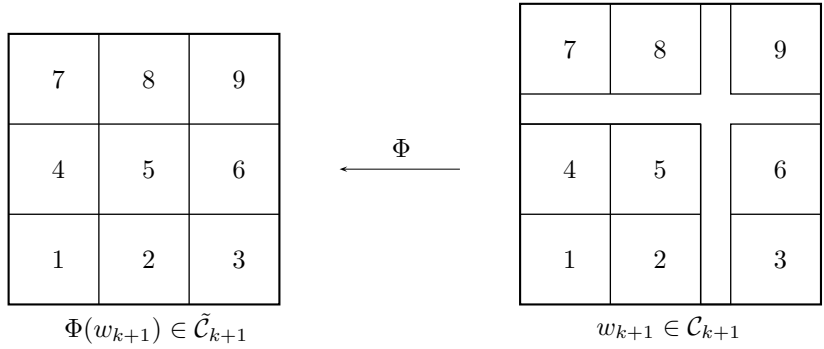
\uput[0](6.5,-2.2){\bf Figure: example of $\Phi$ when $d=2$.}
\end{pspicture}\end{center}

By definition,
$w_{k+1}$ and $\Phi(w_{k+1})$ coincide in $[0,n_{k+1}-n_k-1)^d$, for every $w_{k+1}\in\cc_{k+1}$.
This implies that every element of $\cc_k$ appears in every
element of $\cc_{k+1}$ centered at every residue class modulo $F_k$.

Let us prove that $(X,T)$ is totally minimal. Fix elements $x,y\in X$, a subgroup $F\in\F$ and a pattern
$x|_A$ of $x$ centered in $0\in\Z^d$. By the definition of $X$, $x|_A$ is a subconfiguration of some
$w_k\in\cc_k$, for $k$ large enough such that $F_k<F$. As $y$ is built concatenating elements of $\cc_{k+1}$,
$w_k$ occurs in $y$ centered in every residue class modulo $F$ and the same happens to $x|_A$. In particular,
$x|_A$ occurs in $y$ centered in some $g\in F$, which is exactly the required condition.

\subsection{Total strict ergodicity}
In addition to the ocurrence of configurations in every residue class of subgroups of $\Z^d$, we
also control their frequency. Consider a sequence $(d_k)_{k\ge 1}$ of positive real numbers such that
$\sum_{k\ge1}d_k<+\infty$. Assume that $\cc_1,\ldots,\cc_{k-1},\cc_k$ are non-empty sets satisfying (i),
(ii), (iii), (iv) and
\begin{enumerate}
\item[(vi)] For every $(w_{k-1},w_k,g)\in\cc_{k-1}\times\cc_k\times\Z^d$,
$$\fr(w_{k-1},w_k,g,F_{k-1})\in\left(\dfrac{1-d_{k-1}}{{m_{k-1}}^d\cdot|\cc_{k-1}|}\ \,,\,
\dfrac{1+d_{k-1}}{{m_{k-1}}^d\cdot|\cc_{k-1}|}\right)\cdot$$
\end{enumerate}
Before going to the inductive step, let us make an observation. Condition (vi)
also controls the frequency on subgroups $F$ such that $F_{k-1}<F$. In fact, if $\bar F_{k-1}$
is a complete residue set modulo $F_{k-1}$,
\begin{equation}\label{strict eq 1}
\fr(w_{k-1},w_k,g,F)=\sum_{h\in\bar F_{k-1}\atop{h\equiv_F g}}\fr(w_{k-1},w_k,h,F_{k-1})
\end{equation}
and, as $\left|\{h\in\bar F_{k-1}\,;\,h\equiv_F g\}\right|=(F:F_{k-1})$,
\begin{equation}\label{strict eq 2}
\fr(w_{k-1},w_k,g,F)\in\left(\dfrac{1-d_{k-1}}{(\Z^d:F)\cdot|\cc_{k-1}|}\ \,,\,
\dfrac{1+d_{k-1}}{(\Z^d:F)\cdot|\cc_{k-1}|}\right)\cdot
\end{equation}

We proceed the same way as in the previous subsection: take $l_k$ large and
$\tilde\cc_{k+1}\subset\Omega_{[0,l_km_{k+1})^d}(\cc_k)$ non-empty such that
\begin{equation}\label{strict eq 3}
\fr(w_k,\tilde w_{k+1},g,F_k)\in\left(\dfrac{1-d_k}{{m_k}^d\cdot|\cc_k|}\ \,,\, \dfrac{1+d_k}{{m_k}^d\cdot|\cc_k|}\right)
\end{equation}
for every $(w_k,\tilde w_{k+1},g)\in\cc_k\times\tilde\cc_{k+1}\times\Z^d$.
Note that the non-emptyness of $\tilde\cc_{k+1}$
is guaranteed by Corollary
\ref{law of large numbers}. Also, let $n_{k+1}=l_km_{k+1}n_k+1$ and $\cc_{k+1}=\Phi^{-1}(\tilde\cc_{k+1})$.

Fix $b_A\in\Omega^*(\cc)$. Using the big-$O$ notation, we have
\begin{equation}\label{strict eq 9}
\fr(b_A,W_{k+1},g,F)-\fr(b_A,W_{k+1}|_{[0,n_{k+1}-n_k-1)^d},g,F)=O(1/l_k)\,,
\end{equation}
because these two frequencies differ by the frequency of $b_A$ in $[n_{k+1}-n_k-1,n_{k+1})^d$ and
$$\dfrac{(n_k+1)^d}{{n_{k+1}}^d}=\left(\dfrac{n_k+1}{l_km_{k+1}n_k+1}\right)^d=O(1/l_k).$$
The same happens to $\fr(w_k,w_{k+1},g,F)$ and $\fr(w_k,\Phi(w_{k+1}),g,F)$, because $\Delta(g)=0$
for all $g\in[0,n_{k+1}-n_k-1)^d$. To simplify citation in the future, we write it down:
\begin{equation}\label{strict eq 10}
\fr(w_k,w_{k+1},g,F)-\fr(w_k,\Phi(w_{k+1}),g,F)=O(1/l_k).
\end{equation}
These estimates imply we can assume, taking $l_k$ large enough, that
$$\fr(w_k,w_{k+1},g,F_k)\in\left(\dfrac{1-d_k}{{m_k}^d\cdot|\cc_k|}\ \,,\, \dfrac{1+d_k}{{m_k}^d\cdot|\cc_k|}\right),
\ \forall\,(w_k,w_{k+1},g)\in\cc_k\times\cc_{k+1}\times\Z^d.$$

We make a calculation to be used in the next proposition. Fix $b_A\in\Omega^*(\cc)$ and
$F\in\F$. The main (and simple) observation is: if $b_A$ occurs in $W_k\in\cc_k$
centered at $g$ and $w_k$ occurs in $\Phi(w_{k+1})\in\tilde\cc_{k+1}$ centered at
$h\in[0,l_km_{k+1}-1)^d$, then $b_A$ occurs in $W_{k+1}\in\cc_{k+1}$ centered at $g+n_kh$. This implies
that, if $\bar F$ is a complete residue set modulo $F$, the cardinality of $S(b_A,W_{k+1}|_{[0,n_{k+1}-n_k-1)^d},g,F)$
is equal to
\begin{eqnarray*}
&&\sum_{h\in\bar F\atop{w_k\in\cc_k}}
\sum_{{w_k\text{ occurring in }\atop{{w_{k+1}}|_{[0,l_km_{k+1}-1)^d}\atop{\text{at a vertex }\equiv_F h}}}}
|S(b_A,W_k,g-n_kh,F)|+T\\
&=&\sum_{h\in\bar F\atop{w_k\in\cc_k}}\left|S(w_k,w_{k+1}|_{[0,l_km_{k+1}-1)^d},h,F)\right|
\cdot\left|S(b_A,W_k,g-n_kh,F)\right|+T,
\end{eqnarray*}
where $T$ denotes the number of ocurrences of $b_A$ in $W_{k+1}|_{[0,n_{k+1}-n_k-1)^d}$ not entirely
contained in a concatenated element of $\cc_k$. Observe that\footnote{For each line parallel to a
coordinate axis $e_i\Z$ between two elements of $\cc_k$ in $W_{k+1}$ or containing a line of ones,
there is a rectangle of dimensions $n\times\cdots\times n\times n_{k+1}\times n\times\cdots\times n$
in which $b_A$ is not entirely contained in a concatenated element of $\cc_k$.}
$$0\le T\le d\cdot l_km_{k+1}\cdot n^{d-1}\cdot n_{k+1}<dn^{d-1}\cdot\dfrac{{n_{k+1}}^2}{n_k}\ ,$$
where $n$ is the length of $b_A$. Dividing $\left|S(b_A,W_{k+1}|_{[0,n_{k+1}-n_k-1)^d},g,F)\right|$
by ${n_{k+1}}^d$ and using (\ref{strict eq 9}), (\ref{strict eq 10}), we get
\begin{eqnarray}
\fr(b_A,W_{k+1},g,F)&=&\left(\dfrac{n_{k+1}-1}{n_{k+1}}\right)^d
\sum_{h\in\bar F\atop{w_k\in\cc_k}}\fr(w_k,w_{k+1},h,F)\fr(b_A,W_k,g-n_kh,F)\nonumber\\
&&+\,O(1/l_{k-1})\,.\label{strict eq 4}
\end{eqnarray}
We wish to show that $\fr(b_A,x,F)$ does not depend on $x\in X$. For this, define
\begin{eqnarray*}
\al_k(b_A,F)&=&\min\left\{\fr(b_A,W_k,g,F)\,;\,W_k\in\cc_k,g\in\Z^d\right\}\\
\be_k(b_A,F)&=&\max\left\{\fr(b_A,W_k,g,F)\,;\,W_k\in\cc_k,g\in\Z^d\right\}.
\end{eqnarray*}
The required property is a direct consequence\footnote{In fact, just take the limit in the inequality
$\al_k(b_A,F)\le \fr(b_A,x|_{A_k},0,F)\le\be_k(b_A,F)$.} of the next result.

\begin{proposition}
If $b_A\in\Omega^*(\cc)$ and $F\in\F$, then
$$\lim_{k\rightarrow+\infty}\al_k(b_A,F)=\lim_{k\rightarrow+\infty}\be_k(b_A,F).$$
\end{proposition}

\begin{proof}
By (\ref{strict eq 1}), if $l$ is large such that $F_l<F$, then
$$(F:F_l)\cdot\al_k(b_A,F_l)\le\al_k(b_A,F)\le\be_k(b_A,F)\le(F:F_l)\cdot\be_k(b_A,F_l).$$
This means that we can assume $F=F_l$. We estimate $\al_{k+1}(b_A,F)$ and $\be_{k+1}(b_A,F)$ in
terms of $\al_k(b_A,F)$ and $\be_k(b_A,F)$, for $k\ge l$. As $b_A$ and $F$ are fixed, denote the above
quantities by $\al_k$ and $\be_k$. Take $W_{k+1}\in\cc_{k+1}$.
By (\ref{strict eq 4}),
\begin{eqnarray*}
\fr(b_A,W_{k+1},g,F)&\ge&\left(\dfrac{n_{k+1}-1}{n_{k+1}}\right)^d\cdot\al_k\cdot
\sum_{h\in\bar F\atop{w_k\in\cc_k}}\fr(w_k,w_{k+1},h,F)+O(1/l_{k-1})\\
&=&\left(\dfrac{n_{k+1}-1}{n_{k+1}}\right)^d\cdot\al_k+O(1/l_{k-1})
\end{eqnarray*}
and, as $W_{k+1}$ and $g$ are arbitrary, we get
\begin{equation}\label{strict eq 5}
\al_{k+1}\ge\left(\dfrac{n_{k+1}-1}{n_{k+1}}\right)^d\cdot\al_k+O(1/l_{k-1})\,.
\end{equation}
Equality (\ref{strict eq 4}) also implies the upper bound
\begin{eqnarray}
\fr(b_A,W_{k+1},g,F)&\le&\left(\dfrac{n_{k+1}-1}{n_{k+1}}\right)^d\cdot\be_k\cdot\sum_{h\in\bar F\atop{w_k\in\cc_k}}\fr(w_k,w_{k+1},h,F)+O(1/l_{k-1})\nonumber\\
&=&\left(\dfrac{n_{k+1}-1}{n_{k+1}}\right)^d\cdot\be_k+O(1/l_{k-1})\nonumber\\
\Longrightarrow\hspace{1cm}\be_{k+1}&\le&\left(\dfrac{n_{k+1}-1}{n_{k+1}}\right)^d\cdot\be_k
+O(1/l_{k-1})\,.\label{strict eq 6}
\end{eqnarray}
Inequalities (\ref{strict eq 5}) and (\ref{strict eq 6}) show that $\al_{k+1}$ and $\be_{k+1}$
do not differ very much from $\al_k$ and $\be_k$. The same happens to their difference.
Consider $w_1,w_2\in\cc_{k+1}$ and $g_1,g_2\in\Z^d$. Renaming $g-n_kh$ by $h$ in (\ref{strict eq 4}) and
considering $n_{-k}$ the inverse of $n_k$ modulo $m_l$, the difference
$\fr(b_A,W_1,g_1,F)-\fr(b_A,W_2,g_2,F)$ is at most
\begin{eqnarray*}
&&\sum_{h\in\bar F\atop{w_k\in\cc_k}}\fr(b_A,W_k,h,F)\left|\fr(w_k,w_1,n_{-k}(g_1-h),F)-
\fr(w_k,w_2,n_{-k}(g_2-h),F)\right|\\
&&+\,O(1/l_{k-1})\,.
\end{eqnarray*}
From (\ref{strict eq 3}),
\begin{eqnarray*}
\fr(b_A,W_1,g_1,F)-\fr(b_A,W_2,g_2,F)&\le&\dfrac{2d_k}{{m_l}^d\cdot|\cc_k|}
\sum_{h\in\bar F\atop{w_k\in\cc_k}}\fr(b_A,W_k,h,F)\\
&&+\,O(1/l_{k-1})\\
&\le&2d_k+O(1/l_{k-1})\,,\\
\end{eqnarray*}
implying that
\begin{equation}\label{strict eq 7}
0\ \le\ \be_{k+1}-\al_{k+1}\ \le\ 2d_k+O(1/l_{k-1})\,.
\end{equation}
In particular, $\be_k-\al_k$ converges to zero as $k\rightarrow+\infty$.
The proposition will be proved if $\be_k$ converges. Let us estimate
$|\be_{k+1}-\be_k|$. On one side, (\ref{strict eq 6}) gives
\begin{equation}\label{strict eq 8}
\be_{k+1}-\be_k\ \le\ O(1/l_{k-1})\,.
\end{equation}
On the other, by (\ref{strict eq 5}) and (\ref{strict eq 7}),
\begin{eqnarray*}
\be_{k+1}-\be_k&\ge&\al_{k+1}-\be_k\\
    &\ge&\left(\dfrac{n_{k+1}-1}{n_{k+1}}\right)^d\cdot\al_k-\be_k+O(1/l_{k-1})\\
    &\ge&\left(\dfrac{n_{k+1}-1}{n_{k+1}}\right)^d\cdot\left[\be_k-2d_{k-1}-O(1/l_{k-2})\right]
       -\be_k+O(1/l_{k-1})\nonumber\\
\end{eqnarray*}
which, together with (\ref{strict eq 8}), implies that
\begin{eqnarray*}
\left|\be_{k+1}-\be_k\right|&\le&2d_{k-1}+\be_k\cdot\left[1-\left(\dfrac{n_{k+1}-1}{n_{k+1}}\right)^d\right]+O(1/l_{k-2})\\
    &=&2d_{k-1}+O(1/l_{k-2})\,.
\end{eqnarray*}
As $\sum d_k$ and $\sum 1/l_k$ both converge, $(\be_k)_{k\ge 1}$ is a Cauchy sequence,
which concludes the proof.
\end{proof}

From now on, we consider $(X,T)$ as the dynamical system constructed as above. Note that we have total
freedom to choose $\cc_k$ with few or many elements. This is what controls the entropy of the system.

\subsection{Proof of Theorem \ref{thm 1}}

By (\ref{entropy eq 1}), the topological entropy of the $\Z^d$-action $(X,T)$ satisfies
$$h(X,T)\ge\lim_{k\rightarrow+\infty}\dfrac{\log |\cc_k|}{{n_k}^d}\,\cdot$$
Consider a sequence $(\nu_k)_{k\ge 1}$ of positive real numbers. In the construction of $\cc_{k+1}$
from $\cc_k$, take $l_k$ large enough such that
\begin{enumerate}
\item[(vii)] $\nu_k\cdot n_{k+1}\ge 1$.
\item[(viii)] $|\tilde\cc_{k+1}|\ge |\cc_k|^{(l_km_{k+1})^d\cdot(1-\nu_k)}$.
\end{enumerate}
These inequalities imply
\begin{eqnarray*}
\dfrac{\log |\cc_{k+1}|}{{n_{k+1}}^d}&\ge&\dfrac{\log |\tilde\cc_{k+1}|}{{n_{k+1}}^d}\\
            &\ge&\dfrac{(l_km_{k+1})^d\cdot(1-\nu_k)\cdot\log |\cc_k|}{{n_{k+1}}^d}\\
            &\ge&(1-\nu_k)^{d+1}\cdot\dfrac{\log |\cc_k|}{{n_k}^d}
\end{eqnarray*}
and then
$$\dfrac{\log |\cc_k|}{{n_k}^d}\ge\prod_{i=1}^{k-1}(1-\nu_i)^{d+1}\cdot\dfrac{\log|\cc_1|}{{n_1}^d}
=\prod_{i=1}^{k-1}(1-\nu_i)^{d+1}\cdot\log 2\,.$$
If $\nu\in(0,1)$ is given and $(\nu_k)_{k\ge 1}$ are chosen also satisfying
$$\lim_{k\rightarrow+\infty}\prod_{i=1}^{k}(1-\nu_i)^{d+1}=1-\nu\,,$$
we obtain that $h(X,T)\ge (1-\nu)\log 2>0$. If, instead of $\{0,1\}$, we take $\cc$ with
more elements and apply the construction verifying (i) to (viii), the topological entropy
of the $\Z^d$-action is at least $(1-\nu)\log |\cc|$. We have thus proved Theorem \ref{thm 1}.

\section{Proof of Theorems \ref{thm 3} and \ref{thm 2}}

Given a finite alphabet $\cc$, consider a configuration $b_{A_1}:A_1\rightarrow\cc$ and any $A_2\subset A_1$
such that $|A_2|\le \varepsilon|A_1|$. If $b_{A_2}:A_2\rightarrow\cc$, the element $w\in\Omega_{A_1}(\cc)$
defined by
\begin{eqnarray*}
w(g)&=&b_{A_1}(g)\,,\text{ if }g\in A_1\backslash A_2\,,\\
         &=&b_{A_2}(g)\,,\text{ if }g\in A_2
\end{eqnarray*}
has frequencies not too different from $b_{A_1}$, depending on how small $\varepsilon$ is. In fact,
for any $c\in\cc$,
$$|S(c,b_{A_1},g,F)|-|A_2|\le|S(c,w,g,F)|\le|S(c,b_{A_1},g,F)|+|A_2|$$
and then $|\fr(c,b_{A_1},g,F)-\fr(c,w,g,F)|\le\varepsilon$.

\begin{definition}
The {\it upper-Banach density} of a set $P\subset\Z^d$ is equal to
$$d^*(P)=\limsup_{n_1,\ldots,n_d\rightarrow+\infty}\dfrac{|P\cap [r_1,r_1+n_1)\times\cdots\times[r_d,r_d+n_d)|}
{n_1\cdots n_d}\,\cdot$$
\end{definition}

Consider a set $P\subset\Z^d$ of zero upper-Banach density. We will make $l_k$ grow quickly such that
any pattern of $P\cap(A_k+g)$ appears as a subconfiguration in an element of $\cc_k$. Let's explain
this better. Consider the $d$-dimensional cubes $(A_k)_{k\ge 1}$ that define $(X,T)$. For each $k\ge 1$,
let $\tilde A_k\subset A_k$ be the region containing concatenated elements
of $\cc_{k-1}$. Inductively, they are defined as $\tilde A_1=\{0\}$ and
$$\tilde A_{k+1}=\bigcup_{g\in[0,l_km_{k+1})^d}\left(\tilde A_k+n_kg+\Delta(n_kg)\right)\,,\ \forall\,k\ge 1,$$
where $\Delta$ is the function defined in (\ref{equacao 10}).

\begin{lemma}\label{essential lemma}
If $P\subset\Z^d$ has zero upper-Banach density, there exists a totally strictly ergodic $\Z^d$-action
$(X,T)$ with the following property: for any $k\ge 1$, $g\in\Z^d$ and $b:P\cap(A_k+g)\rightarrow\{0,1\}$,
there exists $w_k\in\cc_k$ such that
$$w_k(h-g)=b(h)\,,\ \forall\,h\in P\cap(A_k+g).$$
\end{lemma}

\begin{proof}
We proceed by induction on $k$. The case $k=1$ is obvious, since $\cc_1\cong\{0,1\}$. Suppose the result
is true for some $k\ge 1$ and consider $b:P\cap(A_{k+1}+g_0)\rightarrow\{0,1\}$. By definition, any $0$,$1$
configuration on $A_{k+1}\backslash\tilde A_{k+1}$ is admissible, so that we only have to worry
about positions belonging to $\tilde A_{k+1}$. For each $g\in[0,l_km_{k+1})^d$, let
$$b^g:P\cap\left(\tilde A_k+n_kg+\Delta(n_kg)+g_0\right)\rightarrow\{0,1\}$$
be the restriction of $b$ to $P\cap\left(\tilde A_k+n_kg+\Delta(n_kg)+g_0\right)$. If $\varepsilon>0$
is given and $l_k$ is large enough,
\begin{eqnarray*}
\dfrac{\left|P\cap\left(\tilde A_{k+1}+g_0\right)\right|}{\left|\tilde A_{k+1}+g_0\right|}&<&
                                          \dfrac{\varepsilon}{(2n_k)^d}\\
\Longrightarrow\hspace{1cm} \left|P\cap\left(\tilde A_{k+1}+g_0\right)\right|&<&\varepsilon\cdot(l_km_{k+1})^d\,,
\end{eqnarray*}
for any $g_0\in\Z^d$. This implies that $P\cap\left(\tilde A_k+n_kg+\Delta(n_kg)+g_0\right)$ is non-empty for at most
$\varepsilon\cdot(l_km_{k+1})^d$ values
of $g\in[0,l_km_{k+1})^d$. For each of these, the inductive hypothesis guarantees the existence of
$w^g\in\cc_k$ such that
$$w^g(h-n_kg-\Delta(n_kg)-g_0)=b^g(h)\,,\ \forall\,h\in P\cap\left(\tilde A_k+n_kg+\Delta(n_kg)+g_0\right).$$
Take any element $z\in\cc_{k+1}$ and define $\tilde z\in\Omega_{A_{k+1}}(\{0,1\})$ by
\begin{eqnarray*}
\tilde z(h)&=&w^g\left(h-n_kg-\Delta(n_kg)\right)\,,\text{ if }h\in \tilde A_k+n_kg+\Delta(n_kg)\\
           &=&b(h)\hspace{3cm},\text{ if }h\in A_{k+1}\backslash\tilde A_{k+1}\\
         &=&z(h)\hspace{3cm},\text{ otherwise}.
\end{eqnarray*}
If $\varepsilon>0$ is sufficiently small, $\tilde z\in\cc_{k+1}$. By its own definition,
$\tilde z$ satisfies the required conditions.
\end{proof}

The above lemma is the main property of our construction. It proves the following stronger statement.

\begin{corollary}\label{main corollary}
Let $(X,T)$ be the $\Z^d$-action obtained by the previous lemma. For any $b:P\rightarrow\{0,1\}$,
there is $x\in X$ such that $x|_P=b$. Also, given $x\in X$, $A\in\mathcal R$ and
$b:P\rightarrow\{0,1\}$, there are $\tilde x\in X$ and $n\in\N$ such that $\tilde x|_A=x|_A$ and
$\tilde x(g)=b(g)$ for all $g\in P\backslash (-n,n)^d$.
\end{corollary}

\begin{proof}
The first statement is a direct consequence of Lemma \ref{essential lemma} and a diagonal argument.
For the second, remember that $x$ is the concatenation of elements of $\cc_k$ and lines of zeroes and ones,
for every $k\ge 1$. Consider $k\ge 1$ sufficiently large and $z_k\in\cc_k$ such that $x|_A$ occurs in $z_k$.
For any $z\in\cc_{k+1}$, there is $g\in\Z^d$ such that $z|_{A_k+g}=z_k$. Constructing $\tilde z$ from $z$
making all substitutions described in Lemma \ref{essential lemma}, except in the pattern $z|_{A_k+g}$, we
still have that $\tilde z\in\cc_{k+1}$.
\end{proof}

\subsection{Proof of Theorem \ref{thm 3}}

Consider $f:X\rightarrow \R$ given by $f(x)=x(0)$. Then
$$\dfrac{1}{\left|P\cap(-n,n)^d\right|}\sum_{g\in P\cap(-n,n)^d}f\left(T^gx\right)=\fr\left(1,x|_{P\cap(-n,n)^d}\right).$$
For each $n\ge 1$, consider the sets
\begin{eqnarray*}
\Lambda_n&=&\bigcup_{k\ge n}\left\{x\in X\,;\,\fr\left(1,x|_{P\cap(-k,k)^d}\right)<1/n\right\}\\
\Lambda^n&=&\bigcup_{k\ge n}\left\{x\in X\,;\,\fr\left(1,x|_{P\cap(-k,k)^d}\right)>1-1/n\right\}.
\end{eqnarray*}
Fixed $k$ and $n$, the sets $\left\{x\in X\,;\,\fr(1,x|_{P\cap(-k,k)^d})<1/n\right\}$ and
is clearly open, so that the same happens to $\Lambda_n$. It is also dense in $X$,
as we will now prove. Fix $x\in X$
and $\varepsilon>0$. Let $k_0\in\N$ be large enough so that
$d(x,y)<\varepsilon$ whenever $x|_{(-k_0,k_0)^d}=y|_{(-k_0,k_0)^d}$. Take $y\in X$ such that
$y|_{(-k_0,k_0)^d}=x|_{(-k_0,k_0)^d}$ and $y(g)=0$ for all $g\in P\backslash(-n,n)^d$ as in
Corollary \ref{main corollary}.
As $\fr(1,y|_{(-k,k)^d})$ approaches to zero as $k$ approaches to infinity,
$y\in\Lambda_n$, proving that $\Lambda_n$ is dense in $X$. The same argument show that $\Lambda^n$
is a dense open set. Then
$$X_0=\bigcap_{n\ge 1}\left(\Lambda_n\cap\Lambda^n\right)$$
is a countable intersection of dense open sets, thus residual. For each $x\in X_0$,
\begin{eqnarray*}
\liminf_{n\rightarrow+\infty}\dfrac{1}{|P\cap(-n,n)^d|}\sum_{g\in P\cap(-n,n)^d}f\left(T^gx\right)&=&0\\
\limsup_{n\rightarrow+\infty}\dfrac{1}{|P\cap(-n,n)^d|}\sum_{g\in P\cap(-n,n)^d}f\left(T^gx\right)&=&1\,,
\end{eqnarray*}
which concludes the proof of Theorem \ref{thm 3}.

\subsection{Proof of Theorem \ref{thm 2}}

Choose an infinite set $G=\{g_i\}_{i\ge 1}$ in $\Z^d$ disjoint from $P$ such that $P'=G\cup P\cup\{0\}$
also has zero upper-Banach density and let $(X,T)$ be the $\Z^d$-action given by Lemma \ref{essential lemma}
with respect to $P'$, that is: for every $b:P'\rightarrow\{0,1\}$, there exists $x^b\in X$ such that
$x^b|_{P'}=b$. Consider
$$X_0=\left\{x^b\in X\,;\,b(0)=0\text{ and }b(g)=1,\ \forall\,g\in P\right\}.$$
This is an uncountable set (it has the same cardinality of $2^G=2^\N$) and, for every $x^b\in X_0$ and $g\in P$,
the elements $T^gx^b$ and $x^b$ differ at $0\in\Z^d$, implying that $x^b\not\in\overline{\{T^gx^b\,;\,g\in P\}}$.
This concludes the proof.

\section*{Acknowledgments}
I would like to thank Vitaly Bergelson for suggesting the topic and for his guidance/advices and enormous
optimism during my visit to The Ohio State University; to The Ohio State University for its great hospitality;
to Carlos Gustavo Moreira, Enrique Pujals and Marcelo Viana for their mathematical support at IMPA;
to Faperj-Brazil for its financial support; to the referee for pointing out some errors.

\bibliographystyle{amsplain}

\end{document}